\documentclass[12pt]{article}
\usepackage{graphicx,psfrag,epsfig}
\usepackage{amssymb,amsmath,amscd,amsthm}
\usepackage{graphicx,psfrag,epsfig}

\usepackage{graphicx}

\newtheorem{theorem}{Theorem}[section]

\newtheorem{lemma}[theorem]{Lemma}

\date{}

\begin{document}

\date{}
\title{A Solvable Model for Homopolymers and the Critical Phenomena}
\author{M. Cranston\footnote{Dept of Mathematics, University of
California, Irvine, CA 92697, mcransto@math.uci.edu}, L.
Koralov\footnote{Dept of Mathematics, University of Maryland,
College Park, MD 20742, koralov@math.umd.edu}, S.
Molchanov\footnote{Dept of Mathematics, University of North
Carolina, Charlotte, NC 28223, smolchan@uncc.edu}, B.
Vainberg\footnote{Dept of Mathematics, University of North
Carolina, Charlotte, NC 28223, brvainbe@uncc.edu}} \maketitle

%
%

\begin{abstract}
We consider a model for the distribution of a long homopolymer
with a zero-range potential at the origin in $\mathbb{R}^3$. The
distribution can be obtained as a limit of Gibbs distributions
corresponding to properly normalized potentials concentrated in
small neighborhoods of the origin as the size of the neighborhoods
tends to zero. The distribution depends on  the length $T$ of the
polymer and a parameter $\gamma$ that corresponds, roughly
speaking, to the difference between the inverse temperature in our
model and the critical value of the inverse temperature.

At the critical point $\gamma_{cr} = 0$ the transition occurs from
the globular phase (positive recurrent behavior of the polymer,
$\gamma > 0$) to the extended phase (Brownian type behavior,
$\gamma < 0$). The main result of the paper is a detailed analysis
of the behavior of the polymer when $\gamma$ is near
$\gamma_{cr}$.

Our approach is based on analyzing the
 semigroups generated by the self-adjoint
extensions $\mathcal{L}_\gamma$ of the Laplacian on $C_0^\infty(
\mathbb{R}^3 \setminus \{ 0\} )$ parametrized by $\gamma$, which
are related to the distribution of the polymer. The main technical
tool of the paper is the explicit formula for the resolvent of the
operator $\mathcal{L}_\gamma$.
\end{abstract}



\section{Introduction}

The study of polymer models has been a very active area of
research in mathematical physics in recent years. Many of the
physically relevant problems have been outlined in the paper of
Lifschitz, Grosberg and Khokhlov \cite{LGK}. To name just a few
recent articles on the subject, with apologies to the many authors
who have been omitted, we cite \cite{AS}, \cite{CGH}, \cite{CH},
\cite{CSY1}, \cite{CSY2}, \cite{CY}. In addition, there is an
interesting exposition with many valuable references in~\cite{G}.
The approach to the subject involves many essential ideas from
statistical mechanics including strong connections to the
developments in \cite{D}, \cite{K} and \cite{LL}.

Let us give a general qualitative description of the problem. We
start with a path space $\Omega = C([0,T], \mathbb{R}^3)$, a
Hamiltonian $H_T:\Omega \rightarrow (-\infty,0]$ and a probability
measure on $\Omega$. A continuous function $\omega \in \Omega$
will be thought of as a realization of the polymer. The parameter
$t \in [0,T]$ can be intuitively understood as the length along
the polymer (although the functions $\omega = \omega(t)$ are not
differentiable and the genuine notion of length can not be
defined). One can consider models with the Hamiltonian given by
\[
H_T(\omega)=-\int_0^T \int_0^T v(\omega(s_1) - \omega(s_2))d s_1
ds_2,
\]
where $v: \mathbb{R}^3 \times \mathbb{R}^3 \rightarrow [0,\infty)$
is a local attractive potential (see \cite{BS}, \cite{KW}). Our
paper concerns a simpler ``mean field" type model, where the
polymer chain interacts with the external attractive potential (as
in \cite{LGK}). Namely, consider the Hamiltonian
\begin{equation} \label{hamt}
H_T(\omega)=-\int_0^T v(\omega(s))d s,
\end{equation}
where $v$ is a
nonnegative compactly supported potential that is not identically
equal to zero. The model considered here is that of a homopolymer
since the potential is time-independent. Other models (of
heteropolymers) which we do not consider in this paper allow for a
potential that is a random stationary function of time  (see
\cite{AS} for example).

Let $\mathrm{P}_{T}^x$ be the Wiener measure on $\Omega = C([0,T],
\mathbb{R}^3)$ shifted so that the trajectories start at $x$
almost surely. This will be the reference measure corresponding to
the infinite temperature, i.e., to the inverse temperature $\beta
= 0$. For a value of inverse temperature $\beta \geq 0$, the
polymer is distributed according to the Gibbs measure
$\mathrm{P}_{\beta,T}^x$, whose density with respect to
$\mathrm{P}_{T}^x$ is
\begin{equation} \label{gibm}
\frac{d \mathrm{P}^x_{\beta,T}}{d \mathrm{P}_{T}^x}(\omega) =
 \frac{\exp(- \beta H_T(\omega))}{%
Z_{\beta}(T,x)},~~\omega \in C([0,T], {\mathbb{R}}^3).
\end{equation}
 The normalizing factor $Z_{\beta}(T,x) = \mathrm{E}_{T}^x
e^{-\beta H_T}$ is called the partition function.

 The phase transitions for the polymers are
well-understood at the physical level (see \cite{LGK}). For large
$\beta$ ($\beta > \beta_{cr}$ for some $\beta_{cr} > 0$) and any
$t \in [0,T]$,  with high probability with respect to the measure
$\mathrm{P}^x_{\beta,T}$ the values of $\omega(t)$ are bounded by
a constant that is independent of $t$. Thus the polymer is in the
globular state. However, for $\beta < \beta_{cr}$ the typical
shape of $\omega(t)$, $t \in [0,T]$, is that of a Brownian path in
$ \mathbb{R}^3$, and the polymer is said to be in the diffusive
state. Rigorous results of this nature have also been proven in
some generality (see \cite{G}, \cite{AS} for example). Of
particular interest is the behavior of the polymer for the values
of $\beta$ that are near $\beta_{cr}$. Unfortunately, the detailed
analysis of critical phenomena for general Hamiltonians is largely
outside the range of techniques of modern mathematical physics.
Most of the results on this type of problems concern the existence
of a phase transition.

However, in our situation (homopolymer with a compactly supported
non-negative potential), the complete analysis of the critical
phenomena is possible. In \cite{CKMV} we studied  the prevalent
behavior of the polymer with respect to the measure
$\mathrm{P}_{\beta,T}^x$ as $T \rightarrow +\infty$. In
particular, we saw that for $d \geq 3$ there is $\beta_{cr} > 0$
such that:

(a) For $\beta > \beta_{cr}$, each constant $s > 0$ and each
function $S(T)$ such that $S(T) \rightarrow +\infty$ and $T - S(T)
\rightarrow +\infty$ as $T \rightarrow +\infty$,  the distribution
of $\omega(S(T)+t)$, $t \in [0, s]$, with respect to
$\mathrm{P}_{\beta,T}^x$, converges to a recurrent Markov process
on $[0,s]$.

(b) For $\beta < \beta_{cr}$, the distribution of $\omega(t
T)/\sqrt{T}$, $t \in [0,1]$, converges to the Brownian motion on
the interval $[0,1]$.

(c) For $\beta = \beta_{cr}$, the distribution of $\omega(t
T)/\sqrt{T}$, $t \in [0,1]$, converges to the distribution of a
non-Gaussian Markov process on the interval $[0,1]$.

Some of the earlier results
for the discrete model based on the random walk and zero-range
potential can be found in \cite{CM}.

In the current paper we study the distributions of the continuous
polymer in $ \mathbb{R}^3$ which correspond to Schrodinger
operators with zero-range potential. The potential is, in some
sense, supported at the origin. Note, however, that the
Hamiltonian $H_T$, the corresponding Gibbs measure and the
Schrodinger operator $L_\beta$ defined by (\ref{hamt}),
(\ref{gibm}) and (\ref{smpl}) are meaningless if $v$ is a
delta-function concentrated at the origin. Before we define these
objects, let us note that in the case of a smooth potential, the
finite-dimensional distributions of $\mathrm{P}^x_{\beta,T}$ can
be expressed in terms of the fundamental solutions
$p_\beta(t,x,y)$ of the parabolic equation
\begin{equation} \label{smpl}
\frac{\partial u}{\partial t} = L_\beta u := \frac{1}{2} \Delta u
+ \beta v u.
\end{equation}
Indeed, by (\ref{gibm}) and the Feynman-Kac formula,
\[
{Z}_{\beta}(T,x) = \int_{ \mathbb{R}^3} {p}_{\beta}(T,x,y) d y,
\]
and
\begin{equation} \label{fkatz}
\mathrm{P}_{\beta,T}^x (\omega(t_1)  \in A_1,...,\omega(t_k) \in
A_k)=
\end{equation}
\[
Z_{\beta}^{-1}( T,x)\int_{A_1}...\int_{A_k} \int_{ \mathbb{R}^3}
p_\beta(t_1,x,x_1)...p_\beta(t_k - t_{k-1}, x_{k-1}, x_k)
p_\beta(T - t_k, x_k, y) d y d x_k...dx_1,
\]
where $k \geq 1$, $0 \leq t_1 \leq ... \leq t_k \leq T$ and
$A_1,...,A_k$ are Borel sets in $ \mathbb{R}^3$. Here we use the
convention that $p_\beta(0,x,y) = \delta_x(y)$. The measure
$\mathrm{P}_{\beta,T}^x$ corresponds to a non-homogeneous Markov
process on $[0,T]$ (see \cite{CKMV}).

Now, instead of the operators $L_\beta$,  we start with the
Schrodinger operators with zero-range potential. They are defined
(see \cite{A}) as self-adjoint extension in $L^2( \mathbb{R}^3)$
of the operator
\[
\Delta/2: C_0^\infty( \mathbb{R}^3 \setminus \{ 0\} ) \rightarrow
C_0^\infty( \mathbb{R}^3 \setminus \{ 0\} ).
\]
There is a family of such extensions $\mathcal{L}_\gamma$ which
depend on a parameter $\gamma \in \mathbb{R}$ that plays the same
role as the difference between  $\beta$ and $\beta_{cr}$ in the
case of smooth potentials (see Section~\ref{pp}). Then we replace
operator $L_\beta$ in (\ref{smpl}) by the operator
$\mathcal{L}_\gamma$. This allows us to define the finite
dimensional distributions and construct the Gibbs measure which
corresponds to the zero-range potentials. Namely, let
$\overline{p}_\gamma(t,x,y)$, $ t \geq 0$, be the kernels of the
operators in the semi-group generated by $\mathcal{L}_\gamma$. We
can now define
\[
\overline{Z}_{\gamma}(T,x) = \int_{ \mathbb{R}^3}
\overline{p}_{\gamma}(T,x,y) d y,
\]
and use an analogue of (\ref{fkatz}) with $\overline{p}_{\gamma}$
instead of $p_\beta$ and $\overline{Z}_{\gamma}(T,x)$ instead of
$Z_{\beta}(T,x)$ in order to define the finite-dimensional
distributions of the measure $ \overline{\mathrm{P}}_{\gamma,
T}^x$. (While neither $ \overline{p}_\gamma(t_1,x,x_1)$ nor
$\overline{Z}_{\gamma}(T,x)$ are defined when $x=0$,  we can make
sense of the expression $ \overline{Z}^{-1}_{\gamma}(T,0)
\overline{p}_\gamma(t_1,0,x_1)$  by taking the limit of
$\overline{Z}^{-1}_{\gamma}(T , x) \overline{p}_\gamma(t_1,x,x_1)$
as $|x| \downarrow 0$.) It is not difficult to check that there is
a measure on $C([0,T], \mathbb{R}^3)$ with such finite-dimensional
distributions, and that this measure defines a non-homogeneous
Markov process on $[0,T]$.

There is another way to obtain the measures $
\overline{\mathrm{P}}^x_{\gamma, T}$ corresponding to zero-range
potentials. One can construct $\overline{\mathrm{P}}^x_{\gamma,
T}$ as the weak limit of measures corresponding to bounded
compactly supported potentials as their supports shrink and the
values of the potentials increase in a particular way. Namely, it
turns out that $ \mathcal{L}_\gamma$, $\gamma \in \mathbb{R}$, can
be obtained as a strong resolvent limit, as $\varepsilon
\downarrow 0$, of Schrodinger operators $
\mathcal{L}_\gamma^\varepsilon = \Delta/2 +  v_\gamma^\varepsilon$
with the potentials (see \cite{A})
\begin{equation}
\label{poten} v_\gamma^\varepsilon = (\frac{\pi^2}{8\varepsilon^2}
+ \frac{\gamma}{\varepsilon}) v(\frac{x}{\varepsilon}),~~||v||_{
L^1( \mathbb{R}^3)}= \frac{4\pi}{3}.
\end{equation}
It will be shown in this paper that the resolvent convergence of
the Schrodinger operators
 implies the convergence of the corresponding Gibbs measures.  Consider the Hamiltonian
$H^{\varepsilon}_{\gamma,T}$ given by (\ref{hamt}) with
$v_{\gamma}^\varepsilon$ instead of $v$, and the Gibbs measure
${{\mathrm{P}}}_{\gamma, T}^{x,\varepsilon}$ given by~(\ref{gibm})
with $\beta = 1$ and $H^{\varepsilon}_{\gamma,T}$ instead of
$H_T$.
%
For each $\gamma \in \mathbb{R}$ and $T > 0$ there are limits
\[
\overline{\mathrm{P}}_{\gamma, T}^x = \lim_{\varepsilon \downarrow
0} {{\mathrm{P}}}_{\gamma, T}^{x,\varepsilon},
\]
understood in the sense of weak convergence of measures on
$C([0,T], \mathbb{R}^3)$ (see Section~\ref{pp}).



Unlike the case of a generic smooth potential,  all the relevant
analytic quantities can be found explicitly in the case of the
zero-range  potential polymer model. This allows us not only to
obtain an analogue of the results (a)-(c) above, but also provide
a detailed analysis of the behavior of the polymer distribution
when $\gamma$ tends to its critical value and, simultaneously, $T
\rightarrow +\infty$.

We recall the main facts about the self-adjoint extensions and
prove the convergence of $ {{\mathrm{P}}}_{\gamma,
T}^{x,\varepsilon}$ to $ {\overline{\mathrm{P}}}_{\gamma,T}^x$ in
Section~\ref{pp}.  We refer the reader to \cite{A} for a detailed
treatment of the self-adjoint extensions of the Laplacian on
$C_0^\infty( \mathbb{R}^3 \setminus \{ 0\} )$. Here we only
mention that for $d \geq 4$ the only closed self-adjoint
extensions of the Laplacian on $C_0^\infty( \mathbb{R}^d \setminus
\{ 0\} )$ is the Laplacian on the entire space, and thus there is
no analog of the distribution corresponding to a zero-range
potential in $d \geq 4$. For $d = 1,2$, the corresponding
homopolymer models do not exhibit a phase transition, so we do not
treat them here.

The goal of the current paper is to analyze the behavior of the
polymer under the measure $\overline{\mathrm{P}}_{\gamma,T}^x$ for
various values of $\gamma$ and large $T$.
Of particular interest is the behavior of the
polymer when $\gamma$ is near zero. We provide a detailed analysis
of the behavior of $\overline{\mathrm{P}}_{\gamma(T),T}^x$ when
$\gamma = \gamma(T) \rightarrow 0$ as $T \rightarrow +\infty$.

Two qualitatively different cases can be considered. The analysis
of both cases relies on a self-similarity property of the measures
$\overline{\mathrm{P}}_{\gamma,T}^x$.  In the first case,
$\gamma(T)$ is bounded and such that  $\gamma(T) \sqrt{T}
\rightarrow +\infty$ as~$T \rightarrow +\infty$. We shall see that
for each constant $s
> 0$ and each function $S(T)$ such that $\gamma(T)\sqrt{S(T)}
\rightarrow +\infty$ and $\gamma(T)\sqrt{T- S(T)} \rightarrow
+\infty$ as $T \rightarrow +\infty$,  the distribution of
$\gamma(T) \omega(S(T)+(\gamma(T))^{-1}t)$, $t \in [0, s]$, with
respect to $\overline{\mathrm{P}}_{\gamma(T),T}^x$, converges to a
recurrent Markov process on $[0,s]$. The generator of the Markov
process will be written out explicitly (see Section~\ref{above}).
In particular, the radial part is a diffusion process on the
positive semi-axis with reflection at the origin. These results
are also applicable in the situation when $\gamma
> 0$ does not depend on $T$.

In the second case $\gamma(T)$ is such that $\gamma(T) \sqrt{T}
\rightarrow \varkappa \in [-\infty, +\infty)$ as~$T \rightarrow
+\infty$. We will show that the distribution of $\omega(t
T)/\sqrt{T}$, $t \in [0,1]$, converges to the limiting measure on
$C([0,1], \mathbb{R}^3)$. The limiting measure can be identified
as $\overline{\mathrm{P}}^{0}_{\varkappa, 1}$ if $\varkappa > -
\infty$ and the Wiener measure on $C([0,1], \mathbb{R}^3)$ if
$\varkappa = -\infty$. In particular, the distribution of
$\omega(t T)/\sqrt{T}$, $t \in [0,1]$, converges to the Wiener
measure on $C([0,1], \mathbb{R}^3)$ if $\gamma < 0$ does not
depend on~$T$.
The limiting distribution for  $\omega(T)/\sqrt{T}$
can be written out explicitly - it turns out to be compound
Gaussian (see Section~\ref{below}).



The paper is organized as follows.
In Section \ref{pp} we recall some facts about the self-adjoint
extensions of the Laplacian and introduce the associated family of
processes with the zero-range potential. In Sections~\ref{above}
and~\ref{below} we  prove the main results concerning the cases
when $\gamma(T) \sqrt{T} \rightarrow +\infty$ and $\gamma(T)
\sqrt{T} \rightarrow \varkappa \in [-\infty, +\infty)$ as~$T
\rightarrow +\infty$, respectively. Finally, in
Section~\ref{tightness} we prove the convergence of the processes
with potentials $v_\gamma^\varepsilon$  to those with a zero-range
potential and prove some technical lemmas regarding the tightness
of certain families of processes.

\section{Measures Corresponding to Zero-Range Potentials} \label{pp}
%

All the self-adjoint extensions of  $\Delta/2$ are described in
the following theorem, whose proof can be found in~\cite{A}.
\begin{theorem}
All the self-adjoint extensions of the Laplacian acting on
$C_0^\infty( \mathbb{R}^3 \setminus \{ 0\} )$ to an operator
acting on $L^2( \mathbb{R}^3)$ form a one-parameter family $
\mathcal{L}_\gamma$, $\gamma \in \mathbb{R}$. The spectrum of $
\mathcal{L}_\gamma$
 is given by
\[
{\rm spec}(
\mathcal{L}_\gamma)=(-\infty,0]\cup\left\{\frac{\gamma^2}{2}\right\},~~\gamma>0,
\]
\[
{\rm spec}( \mathcal{L}_\gamma )=(-\infty,0],~~\gamma \leq 0.
\]
The kernel of the resolvent of $\mathcal{L}_\gamma$ is given by
\[
R_{\lambda,\gamma}(x,y)=
\frac{e^{-\sqrt{2\lambda}|x-y|}}{\pi|x-y|}+\frac{1}{\sqrt{2\lambda}-\gamma}
\frac{e^{-\sqrt{2\lambda}(|x|+|y|)}}{2\pi|x||y|},~~\lambda \notin
{\rm spec}( \mathcal{L}_\gamma).
\]
If $\gamma > 0$, then $\gamma^2/2$ is a simple eigenvalue of
$\mathcal{L}_\gamma$ with the eigenfunction
\[
\psi_{\gamma}(x)=\frac{\sqrt{\gamma}}{\sqrt{2\pi}}\frac
{e^{-\gamma |x|}}{|x|}.
\]
\end{theorem}
Since the spectrum of $ \mathcal{L}_\gamma$ is bounded from above,
the operators $\exp(t \mathcal{L}_\gamma)$, $t \geq 0$,  are
bounded in $L^2( \mathbb{R}^3)$. The kernel of  $\exp(t
\mathcal{L}_\gamma)$, $ t > 0$,  is given by
\begin{equation} \label{funs}
\overline{p}_{\gamma}(t,x,y) = \frac{1}{2\pi i}\int_{\Gamma(a)}
e^{\lambda t} R_{\lambda,\gamma}(x,y) d \lambda =
\frac{e^{-|x-y|^2/2t}}{(2\pi t)^{3/2}}
+\frac{1}{4\pi^2 i}\int_{\Gamma(a)}
\frac{e^{-\sqrt{2\lambda}(|x|+|y|)+\lambda t}}{ (\sqrt{2\lambda}-\gamma)|x||y|} d \lambda,
\end{equation}
where $x, y \neq 0$, $a > \gamma^2/2$ and $\Gamma(a)$ is the
contour in the complex plane that is parallel to the imaginary
axis and passes through $a$.  Thus $\overline{p}_{\gamma}(t,x,y)$
can be interpreted as a formal fundamental solution of the
equation
\[
\frac{\partial u}{\partial t} = \mathcal{L}_\gamma u.
\]
Let
\[
\overline{Z}_{\gamma}(t,x) = \int_{ \mathbb{R}^3}
\overline{p}_{\gamma}(t,x,y) d y,
\]
where $t > 0$, $x \neq 0$. Formula (\ref{funs}) implies that
\begin{equation} \label{zbar}
\overline{Z}_{\gamma}(t , x)=1+\frac{1}{2\pi i}\int_{\Gamma(a)}
e^{\lambda t} \frac{1}{\sqrt{2\lambda}-\gamma}
\frac{e^{-\sqrt{2\lambda}|x|}}{\lambda|x|} d\lambda.
\end{equation}

Define the measures
$\overline{\mathrm{P}}^x_{\gamma, T}$, $x \in \mathbb{R}^3$,  via
their finite-dimensional distributions
\begin{equation} \label{fkatz2}
\overline{\mathrm{P}}^x_{\gamma, T}(\omega(t_1)  \in
A_1,...,\omega(t_k) \in A_k)=
\end{equation}
\[
\overline{Z}_{\gamma}^{-1} (T,x)\int_{A_1}...\int_{A_k} \int_{
\mathbb{R}^3}
\overline{p}_\gamma(t_1,x,x_1)...\overline{p}_\gamma(t_k -
t_{k-1}, x_{k-1}, x_k) \overline{p}_\gamma(T - t_k, x_k, y) d y d
x_k...dx_1,
\]
where $k \geq 1$, $0 \leq t_1 \leq ... \leq t_k \leq T$ and
$A_1,...,A_k$ are Borel sets in $ \mathbb{R}^3$. Here we use the
conventions that $\overline{p}_\gamma(0,x,y) = \delta_x(y)$ and $
\overline{Z}^{-1}_{\gamma}(T, 0) \overline{p}_\gamma(t_1,0,x_1) =
\lim_{|x| \downarrow 0} \overline{Z}^{-1}_{\gamma}(T, x)
\overline{p}_\gamma(t_1,x,x_1)$.
It is not difficult to show that there indeed  is a measure on
$C([0,T], \mathbb{R}^3)$ with the finite-dimensional distributions
given by (\ref{fkatz2}). We don't prove this here since the same
conclusion follows from the proof of Theorem~\ref{convlem} below.
Note the following self-similarity property.
\begin{theorem}
 For $a > 0$,
let $f_a: C([0,T], \mathbb{R}^3) \rightarrow  C([0,a T],
\mathbb{R}^3)$ map a function $\omega$ to the function $f_a
\omega$ via $(f_a \omega)(t) = \sqrt{a}\omega (t/a)$. Let $f_a
\overline{\mathrm{P}}^x_{\gamma, T}$  be the push-forward of the
measure $\overline{\mathrm{P}}^x_{\gamma, T}$ by this mapping.
Then
\begin{equation} \label{ssim}
\overline{p}_{\gamma}(t a,x\sqrt{a},y\sqrt{a}) = a^{-3/2}
\overline{p}_{\gamma \sqrt{a}}(t,x,y),~~~~~\overline{Z}_{\gamma}
(t a ,x\sqrt{a}) = \overline{Z}_{\gamma \sqrt{a}} (t,x)
\end{equation}
and
\begin{equation} \label{ssim2}
f_a \overline{\mathrm{P}}^x_{\gamma, T} =
\overline{\mathrm{P}}^{x\sqrt{a}}_{\gamma /\sqrt{a}, T a}.
\end{equation}
\end{theorem}
\proof The first statement follows after making the change of
variables $\lambda' = \lambda/a$ in the integrals in the right
hand sides of (\ref{funs}) and (\ref{zbar}). The second statement
now follows from~(\ref{fkatz2}). \qed

Let us discuss an alternative construction of the measure
${{\mathrm{P}}}_{\gamma, T}^{x,\varepsilon}$ that uses short-range
potentials.
Let $v_\gamma^\varepsilon $ be defined by (\ref{poten}). Let the
Hamiltonian $H^{\varepsilon}_{\gamma,T}$ be given by~(\ref{hamt})
with $v_\gamma^\varepsilon$ instead of $v$, and the Gibbs measure
${{\mathrm{P}}}_{\gamma, T}^{x,\varepsilon}$ be given by
(\ref{gibm}) with $\beta = 1$ and $H^{\varepsilon}_{\gamma, T}$
instead of~$H_T$. The following theorem will be proved in
Section~\ref{tightness}.

\begin{theorem} \label{convlem}
Let $x \in \mathbb{R}^3 \setminus 0$, $\gamma \in \mathbb{R}$ and
$T > 0$. The measures ${{\mathrm{P}}}_{\gamma, T}^{x,\varepsilon}$
converge, in the sense of weak convergence of measures on
$C([0,T], \mathbb{R}^3)$, as $\varepsilon \downarrow 0$, to the
measure $ \overline{\mathrm{P}}^x_{\gamma, T}$.
\end{theorem}

From~(\ref{fkatz2}) it follows that the finite-dimensional
distributions of $\overline{\mathrm{P}}^{x'}_{\gamma', T}$
converge to those of $\overline{\mathrm{P}}^0_{\gamma, T}$ as $x'
\rightarrow 0$ and  $\gamma' \rightarrow \gamma \in \mathbb{R}$.
If $\gamma' \rightarrow -\infty$, then they converge to the
finite-dimensional distributions of the Wiener measure on
$C([0,T], \mathbb{R}^3)$. From  Theorem~\ref{convlem} and
Lemma~\ref{lemtight} (from Section~\ref{tightness}) it follows
that the family $\{ \overline{\mathrm{P}}^{x'}_{\gamma', T},
\gamma' \leq c, |x'| \leq 1 \}$ is tight for each $c \in
\mathbb{R}$, and therefore
\begin{equation} \label{cnv}
\lim_{x' \rightarrow 0, \gamma' \rightarrow \gamma}
\overline{\mathrm{P}}^{x'}_{\gamma', T} =
\overline{\mathrm{P}}^{0}_{\gamma, T},~~~\lim_{x' \rightarrow 0,
\gamma' \rightarrow -\infty} \overline{\mathrm{P}}^{x'}_{\gamma',
T} = \mathrm{P}_{T}^0,
\end{equation}
where $ \mathrm{P}_{T}^0$ is the Wiener measure and the limits are
understood in the sense of weak convergence of measures.



\section{Distribution Above the Critical Point} \label{above}
In this section we assume that $\gamma  = \gamma(T)$ is bounded
and such that $\gamma(T) \sqrt{T} \rightarrow +\infty$ as~$T
\rightarrow +\infty$. In particular, this covers the case when
$\gamma > 0$ does not depend on $T$. We examine the behavior of
the polymer paths with respect to
$\overline{\mathrm{P}}^x_{\gamma(T), T}$ as $T \rightarrow
+\infty$. First, we need the asymptotics of $
\overline{p}_{1}(t,x,y)$ and $ \overline{Z}_{1} (t,x)$ when $t
\rightarrow +\infty$.
\begin{lemma} \label{bnbn1}
We have the following asymptotic expressions:
\begin{equation} \label{p1}
\overline{p}_{1}(t,x,y) = \exp(t/2)(\psi_1(x)\psi_1(y) +
q(t,x,y)),
\end{equation} \label{z1}
where $\lim_{t \rightarrow +\infty} \sup_{|x|,|y| \geq k}
|q(t,x,y)| = 0$ for each $k > 0$;
\begin{equation}
\overline{Z}_{1} (t,x) = \exp(t/2)||\psi_1||_{L^1(
\mathbb{R}^3)}(\psi_1(x)+ Q(t,x)),
\end{equation}
where $\lim_{t \rightarrow +\infty} \sup_{|x| \geq k} |Q(t,x)| =
0$ for each $k > 0$;
\begin{equation}
 \overline{Z}^{-1}_{1}
(s,x) \overline{p}_{1}(t,x,y)  = \exp((t-s)/2)||\psi_1||_{L^1(
\mathbb{R}^3)}^{-1}(\psi_1(y) + \widetilde{q}(s,t,x,y)),
\end{equation}
where $\lim_{s,t \rightarrow +\infty} \sup_{|x| \leq k^{-1}, |y|
\geq k} |\widetilde{q}(s,t,x,y)| = 0$ for each $k > 0$.
\end{lemma}

\proof Formula (\ref{p1}) follows immediately from (\ref{funs}) if
the integral over $\Gamma(a)$ in (\ref{funs}) is replaced by the
integral over $\Gamma(b),~0<b<\gamma^2/2$, plus the residue at
$\lambda=\gamma^2/2=1/2.$ Formula (\ref{p1}) follows similarly
from (\ref{zbar}).
 \qed
\\

Next, let us study the distribution of the end of the polymer with
respect to the measure $\overline{\mathrm{P}}^x_{\gamma(T), T}$ as
$T \rightarrow +\infty$.
\begin{theorem} If $\gamma(T)$ is bounded and such that
$\gamma(T) \sqrt{T} \rightarrow +\infty$ as $T \rightarrow
+\infty$, then the distribution of $\gamma(T)\omega(T)$ with
respect to the measure $\overline{\mathrm{P}}^x_{\gamma(T), T}$
converges, weakly, as $T \rightarrow +\infty$, to the distribution
with the density $ \psi_1/||\psi_1||_{L^1( \mathbb{R}^3)}$.
\end{theorem}
\proof  For fixed $x$, the density of $\omega(T)$ with respect to
the Lebesgue measure is equal to
\[
\overline{Z}^{-1}_{\gamma(T)} (T,x)
\overline{p}_{\gamma(T)}(T,x,y).
\]
Therefore, the density of $\gamma(T)\omega(T)$ is equal to
\[
\gamma(T)^{-3}\overline{Z}^{-1}_{\gamma(T)} (T,x)
\overline{p}_{\gamma(T)}(T,x,y/\gamma(T)) = \overline{Z}^{-1}_{1}
(\gamma^2(T) T, x \gamma(T))
\overline{p}_{1}(\gamma^2(T)T,x\gamma(T),y).
\]
Here we used (\ref{ssim}) with $\gamma = 1$, $a = \gamma^2(T)$ and
$t = T$. The latter expression is equal to
\begin{equation} \label{hh01}
||\psi_1||^{-1}_{L^1( \mathbb{R}^3)} (\psi_1(y) +
\widetilde{q}(\gamma^2(T)T,\gamma^2(T)T,x \gamma(T),y)),
\end{equation}
 where  $\widetilde{q}$ is the same as
in Lemma~\ref{bnbn1}. When $T \rightarrow +\infty$, the expression
in (\ref{hh01}) converges to $\psi_1(x)/||\psi_1||_{L^1(
\mathbb{R}^3)}$ uniformly in $|y| > k$ by Lemma~\ref{bnbn1}. This
justifies the weak convergence. \qed

Now let us examine the behavior of the polymer in a region
separated both from zero and $T$. Let $S(T)$ be such that
\begin{equation} \label{sttt}
\lim_{T \rightarrow +\infty} \gamma(T) \sqrt{S(T)} = \lim_{T
\rightarrow +\infty} \gamma(T) \sqrt{T - S(T)} = +\infty.
\end{equation}
Let $s > 0$ be fixed. Consider the process $y^T(t) = \gamma(T)
\omega(S(T)+t/\gamma^2(T))$, $0 \leq t \leq s$. Let
\begin{equation} \label{rprp}
 r(t,x,y) = \frac{\overline{p}_1(t,x,y)
\psi_1(y)}{\psi_1(x)} \exp(-t/2), ~~y \neq 0.
\end{equation}
Observe that $\int_{ \mathbb{R}^3}  r(t,x,y) d y = 1$, and
therefore $ r$ serves as the transition density for a Markov
process. Also note that $\psi_1^2$ serves as an invariant density
for the process.
\begin{theorem} \label{pro1}  If $\gamma(T)$ is bounded and such that
$\gamma(T) \sqrt{T} \rightarrow +\infty$ as $T \rightarrow
+\infty$, then the distribution of the process $y^T(t)$ with
respect to the measure  $\overline{\mathrm{P}}^x_{\gamma(T), T}$
converges as $T \rightarrow +\infty$, weakly in the space
$C([0,s], \mathbb{R}^3)$, to the distribution of a stationary
Markov process with transition density $ r(t,x,y)$ and the
invariant measure whose density is $\psi_1^2$.
\end{theorem}
\proof
 Let $0 \leq t_1 < ... < t_n \leq s$. The density of the random
vector $(y^T(t_1),...,y^T(t_n))$ with respect to the Lebesgue
measure on $ \mathbb{R}^{3n}$ is equal to
\[
\rho^T(x_1,...,x_n) =
\]
\[
\gamma(T)^{-3n}(\overline{Z}_{\gamma(T)}(T,x))^{-1}
\overline{p}_{\gamma(T)}(S(T) + \frac{t_1}{\gamma^2(T)}, x,
\frac{x_1}{\gamma(T)}) \overline{p}_{\gamma(T)}(\frac{t_2-
t_1}{\gamma^2(T)}, \frac{x_1}{\gamma(T)},
\frac{x_2}{\gamma(T)})~...
\]
\[
...~\overline{p}_{\gamma(T)}(\frac{t_n -
t_{n-1}}{\gamma^2(T)},\frac{x_{n-1}}{\gamma(T)},
\frac{x_{n}}{\gamma(T)}) \overline{Z}_{\gamma(T)}(T -S(T) -
\frac{t_n}{\gamma^2(T)}, \frac{x_n}{\gamma(T)}).
\]
Applying (\ref{ssim}) with $\gamma = 1$, $a = \gamma^2(T)$ and $t
= T$, we see that the right hand side of this equality is equal to
\[
(\overline{Z}_{1}(\gamma^2(T) T, x\gamma(T)))^{-1}
\overline{p}_{1}(\gamma^2(T) S(T) + {t_1}, x \gamma(T), {x_1})
\overline{p}_{1}({t_2- t_1},{x_1},{x_2})~...
\]
\[
...~\overline{p}_{1}({t_n - t_{n-1}},{x_{n-1}}, {x_{n}})
\overline{Z}_{1}(\gamma^2(T)( T -S(T)) - t_n, {x_n}).
\]
We replace here all factors $\overline{p}_{1}$, except the first
one, by $r$ using (\ref{rprp}). We replace the factors
$(\overline{Z}_{1}(\gamma^2(T) T, x\gamma(T)))^{-1}
\overline{p}_{1}(\gamma^2(T) S(T) + {t_1}, x \gamma(T), {x_1})$
and $\overline{Z}_{1}(\gamma^2(T)( T -S(T)) - t_n, {x_n})$ by
their asymptotic expansions given in Lemma~\ref{bnbn1}. This leads
to
\[
\rho^T(x_1,...,x_n) = \psi^2_1(x_1) r(t_2- t_1, x_1, x_2)... r(t_n
- t_{n-1},x_{n-1}, x_{n})+o(1),~~~T\rightarrow +\infty,
\]
where the remainder tends to zero uniformly in $(x_1,...,x_n)$
with $\min(|x_1|,...,|x_n|) \geq k$. Since $k>0$ can be chosen to
be arbitrarily small, this justifies the convergence of the finite
dimensional distributions of $y^T$ to those of the Markov process.
It remains to note that the family of measures induced by the
processes $y^T$  is tight, as follows from Lemma~\ref{lemtight2}.
\qed
\\
\\
{\bf Remark.} Let us describe the generator of the limiting Markov
process. Namely, let $C_0( \mathbb{R}^3)$ be the space of
continuous functions on $ \mathbb{R}^3$ with a finite limit at
infinity. Consider the differential operator
\[
M u(x) =  \frac{1}{2}\Delta u(x)  + \frac{(\nabla \psi_1 (x),
\nabla u (x))}{\psi_1(x)} = \frac{1}{2}\Delta u(x)  - (1 +
\frac{1}{|x|}) \frac{\partial u(x)}{\partial r},~~x \neq 0,
\]
\[
M u(0) = \lim_{|x| \downarrow 0} ( \frac{1}{2}\Delta u(x)  - (1 +
\frac{1}{|x|}) \frac{\partial u(x)}{\partial r}).
\]
Note that in the spherical coordinates $M$ can be written as
\[
M u = \frac{1}{2}( \frac{\partial^2 u}{ \partial r^2} +
\frac{1}{r^2}\widetilde{\Delta} u)  - \frac{\partial u}{\partial
r},
\]
where $ \widetilde{\Delta}$ is the Beltrami-Laplace operator on
the sphere.
 For fixed $x$, the
function $ r(t,x,y)$ satisfies
\[
\frac{\partial r(t,x,y)}{\partial t} = M^* r(t,x,y),~~t > 0, ~y
\in \mathbb{R}^3 \setminus 0,
\]
where $M^*$ is the formal adjoint of $M$. From here it easily
follows (\cite{Ko}) that the generator in the space $C_0(
\mathbb{R}^3)$ of the Markov family with transition density
$r(t,x,y)$ is given by the operator $M$ with the domain
\[
\mathcal{D} = \{u \in C^2(\mathbb{R}^3 \setminus \{0\}) \cap
C_0(\mathbb{R}^3):~ M u \in C_0(\mathbb{R}^3) ,~ \lim_{r
\downarrow 0} { \frac{ \int_{S^{2}} u (r, \varphi) d\mu(\varphi) -
u(0)}{r} = 0 } \},
\]
where $\mu$ is the Haar probability measure on $S^{2}$. The radial
part of the limiting process is then a diffusion on the positive
semi-axis with unit diffusion coefficient, unit drift towards the
origin and reflection at the origin. If we denote the radial part
of the limiting process by $R_t$, then the spherical part of the
limiting process $S_t$ on $S^2$ satisfies
\begin{equation} \label{enn}
dS_t = \frac{1}{R_t} dB_t,~~ R_t \neq 0.
\end{equation}
 Here $B_t$ is the diffusion
on the sphere whose generator is one half of the Beltrami-Laplace
operator. It is not difficult to see that if we denote the the
first time when $R_t = 0$ by $\tau$, then the set of limit points
of $S_t$ as $t \uparrow \tau$ coincides with the entire sphere
$S^2$. For $t > 0$ the distribution of $S_{\tau+t}$ is uniform on
the sphere.
\\
\\
{\bf Remark.} If instead of assuming that $\lim_{T \rightarrow
+\infty} \gamma(T) \sqrt{S(T)} = +\infty$  we assume that $S(T) =
0$, the result of Theorem~\ref{pro1} will hold with the only
difference that the initial distribution for the limiting Markov
process will now be concentrated at $ x \lim_{T \rightarrow
+\infty} \gamma(T) $, instead of being the invariant distribution,
provided that the latter limit exists.

\section{Distribution Near and Below the Critical Point}
\label{below} In this section we assume that $\gamma  = \gamma(T)$
is such that $\gamma(T) \sqrt{T} \rightarrow \varkappa \in [-\infty,
+\infty)$ as~$T \rightarrow +\infty$. In particular, this covers
the case when $\gamma < 0$ does not depend on $T$. We examine the
behavior of the polymer paths with respect to
$\overline{\mathrm{P}}^x_{\gamma(T), T}$ as $T \rightarrow
+\infty$.

Consider the process $y^T(t) = \omega(t T)/\sqrt{T}$, $0 \leq t
\leq 1$.
Applying (\ref{ssim2}) with $a = 1/T$, we see that
the  distribution of the process $y^T$ with respect to the measure
$\overline{\mathrm{P}}^x_{\gamma(T), T}$ coincides with the
distribution of the process $\omega(t)$, $t \in [0,1]$, with
respect to the
measure~$\overline{\mathrm{P}}^{x/\sqrt{T}}_{\gamma(T)
\sqrt{T},1}$. Since $x/\sqrt{T} \rightarrow 0$ and $\gamma(T)
\sqrt{T} \rightarrow \varkappa$ as $T \rightarrow +\infty$, the following
theorem is a consequence of (\ref{cnv}).
%
\begin{theorem}
If  $ \gamma = \gamma(T)$ is such that  $\gamma(T) \sqrt{T}
\rightarrow \varkappa \in (-\infty, +\infty)$ as $T \rightarrow +\infty$,
then the distribution of the process $y^T(t)$ with respect to the
measure $\overline{\mathrm{P}}^x_{\gamma(T), T}$ converges as $T
\rightarrow +\infty$ to the measure $\overline{\mathrm{P}}^{0}_{\varkappa,
1}$.

If   $ \gamma = \gamma(T)$ is such that  $\gamma(T) \sqrt{T}
\rightarrow -\infty$ as $T \rightarrow +\infty$, then the
distribution of the process $y^T(t)$ with respect to the measure
$\overline{\mathrm{P}}^x_{\gamma, T}$ converges as $T \rightarrow
+\infty$ to the distribution of the $3$-dimensional Brownian
motion.
\end{theorem}
The finite-dimensional distributions of the limiting process can
be written out using (\ref{funs}) and (\ref{fkatz2}). In
particular, when $\varkappa = 0$, the distribution of the limiting
process is the same as the distribution of the corresponding
process in the case of the smooth potential with $\beta =
\beta_{cr}$ (see \cite{CKMV}).

Let us show that the limiting distribution of $y^T(T)$ is compound
Gaussian.
\begin{theorem}
The distribution of $\omega(1)$ (induced by the measure
$\overline{\mathrm{P}}^{0}_{\varkappa,1}$) is compound Gaussian,
i.e., its density is given by
\[
q_{\varkappa }(y)=\int_{0}^{1}\frac{ e^{-|y|^{2}/2\tau }}{(2\pi
\tau )^{3/2}}v_{\varkappa}(\tau )d\tau,
\]
where $v_{\varkappa}(\tau ) > 0$, $\tau \in [0,1]$, $\varkappa \in
\mathbb{R}$.
\end{theorem}
\proof By (\ref{cnv}), (\ref{zbar}) and (\ref{fkatz2}), the
density of $\omega(1)$ with respect to the Lebesgue measure is
equal to
\[
q_{\varkappa }(y) =
\lim_{|x| \downarrow 0}
\frac{ \overline{p}_{\varkappa }(1,x,y) }{\overline{Z}_{\varkappa}(1, x)}=
\frac{1}{2\pi}\int_{\Gamma (a)}\frac{e^{-\sqrt{2\lambda }
|y|+\lambda }}{(\sqrt{2\lambda }-\varkappa )|y|}d\lambda
 \left( \int_{\Gamma(a)}
 \frac{e^{\lambda}d\lambda}{(\sqrt{2\lambda}-\varkappa)\lambda}
 \right)^{-1},
\]
where $a>\varkappa
^{2}/2$. Note that
\begin{equation*}
\frac{e^{-\sqrt{2\lambda }|y|}}{2\pi |y|}=\int_{0}^{\infty }\frac{%
e^{-|y|^{2}/2\tau }}{(2\pi \tau )^{3/2}}e^{-\lambda \tau }d\tau
=\int_{0}^{\infty }\frac{d}{d\tau}\left( \frac{e^{-|y|^{2}/2\tau }}{(2\pi \tau )^{3/2}}%
\right)\frac{e^{-\lambda \tau }}{\lambda }d\tau ,~\lambda \in \Gamma(a).
\end{equation*}
We put this expression in the first integral
in the previous formula. The integration by parts above was needed
to guarantee the absolute convergence of the double integral and to change
the order of integration. Below we will move the derivative
back to the second factor. Thus
\begin{equation} \label{nmnm1}
q_{\varkappa }(y)=\int_{0}^{\infty }\frac{d}{d\tau }\left( \frac{%
e^{-|y|^{2}/2\tau }}{(2\pi \tau )^{3/2}}\right) u_{\varkappa}(\tau
)d\tau ,\text{ \ \ \ \ }u_{\varkappa}(\tau )= K \int_{\Gamma
(a)}\frac{e^{\lambda (1-\tau )}d\lambda }{(\sqrt{2\lambda
}-\varkappa )\lambda },
\end{equation}
where
\[
K =
\left( \int_{\Gamma(a)}
 \frac{e^{\lambda}d\lambda}{(\sqrt{2\lambda}-\varkappa)\lambda}
 \right)^{-1}.
\]
If $\tau \geq 1$ and Re$\lambda \geq a,$ then the integrand in the
second integral in (\ref{nmnm1}) is analytic in $\lambda $ and of
order $O(|\lambda |^{-3/2})$ as $|\lambda |\rightarrow \infty .$
Thus $u_{\varkappa}(\tau )=0$ for $\tau \geq 1,$ and
\begin{equation}
q_{\varkappa }(y)=\int_{0}^{1}\frac{d}{d\tau }\left( \frac{e^{-|y|^{2}/2\tau }}{%
(2\pi \tau )^{3/2}}\right) u_{\varkappa}(\tau )d\tau =\int_{0}^{1}\frac{%
e^{-|y|^{2}/2\tau }}{(2\pi \tau )^{3/2}}v_{\varkappa}(\tau )d\tau,
\label{sos}
\end{equation}
where
\begin{equation*}
v_{\varkappa}(\tau )=-u_{\varkappa}^{\prime }(\tau )=K \int_{\Gamma (a)}\frac{%
e^{\lambda (1-\tau )}d\lambda }{(\sqrt{2\lambda }-\varkappa )},\text{ \ \ \ }%
0<\tau <1.
\end{equation*}
Since the integral over $\Gamma (a)$ is equal to the integral over a contour
around negative $\lambda $-semiaxis (plus the contribution from the residue at $\varkappa ^{2}/2$ if $\varkappa > 0$), we
have
\begin{equation}
\int_{\Gamma (a)}\frac{%
e^{\lambda (1-\tau )}d\lambda }{(\sqrt{2\lambda }-\varkappa )}=2i\int_{0}^{\infty
}\frac{\sqrt{2\sigma }e^{-\sigma (1-\tau )}d\sigma }{2\sigma
+\varkappa ^{2}}+\delta 2\pi i\varkappa e^{\varkappa ^{2}(1-\tau )/2},\text{ \ \ \
}0<\tau <1,  \label{sos1}
\end{equation}
where $\delta =1$ if $\varkappa >0$, $\delta =0$ if $\varkappa \leq 0$.

Let us evaluate $K$. We have
\begin{equation} \label{k-1}
K^{-1}=\int_{\Gamma(a)}
 \frac{e^{\lambda}d\lambda}{(\sqrt{2\lambda}-\varkappa)\lambda}=\int_{\Gamma(a)}
 \frac{(\sqrt{2\lambda}+\varkappa)e^{\lambda}d\lambda}{({2\lambda}-\varkappa ^{2})\lambda}=
\int_{\Gamma(a)}
 \frac{\sqrt{2\lambda}e^{\lambda}d\lambda}{({2\lambda}-\varkappa ^{2})\lambda}+
 \int_{\Gamma(a)}
 \frac{\varkappa e^{\lambda}d\lambda}{({2\lambda}-\varkappa ^{2})\lambda}.
\end{equation}
The first integral $I_1$ in the right-hand side can be written as
the integral over the contour around the negative $\lambda
$-semiaxis plus the contribution from the residue at $\varkappa
^{2}/2$, i.e.
\[
I_1=-2\sqrt{2}i \int_{0}^{\infty}
 \frac{ e^{-\sigma}d\sigma}{({2\sigma}+\varkappa ^{2})\sqrt{\sigma}}+2\pi i\frac{e^{\varkappa^2/2}}{|\varkappa|}.
\]
We make the substitution $\sigma\rightarrow \sigma^2$ in the first
 term of the right-hand side and rewrite the second term using the identity
\[
\frac{\pi}{|\varkappa|}= \int_{0}^{\infty}\frac{\sqrt{2} d\sigma}{\sigma^2+\varkappa ^{2}/2}.
\]
The second integral in the right-hand side of (\ref{k-1}) can be expressed
through the residues at $\lambda =0$ and $\lambda =\varkappa ^{2}/2$. This implies, that
\[
K^{-1}=2\sqrt{2}i\int_{0}^{\infty}\frac{e^{\varkappa^2/2}-e^{-\sigma^2}}{\sigma^2+\varkappa ^{2}/2}d\sigma+
2\pi i\frac{e^{\varkappa^2/2}-1}{\varkappa}.
\]
This and (\ref{sos1}) show that $v_{\varkappa}(\tau )>0$ for all
$0\leq \tau \leq 1,~ -\infty < \varkappa < \infty$. This together
with (\ref{sos}) justify the statement. \qed

\section{Proof of Theorem~\ref{convlem}. Tightness of Certain Families of Processes}
\label{tightness} {\it Proof of Theorem~\ref{convlem}.} As follows
from (\ref{smpl}) and (\ref{fkatz}), the finite-dimensional
distributions of ${{\mathrm{P}}}_{\gamma, T}^{x,\varepsilon}$ are
given by
\[
{{\mathrm{P}}}_{\gamma, T}^{x,\varepsilon}(\omega(t_1)  \in
A_1,...,\omega(t_k) \in A_k)=
\]
\[
({\widetilde{Z}}_{\gamma}^{\varepsilon})^{-1}
(T,x)\int_{A_1}...\int_{A_k} \int_{ \mathbb{R}^3}
{\widetilde{p}}^\varepsilon_\gamma(t_1,x,x_1)...{\widetilde{p}}^\varepsilon_\gamma(t_k
- t_{k-1}, x_{k-1}, x_k) {\widetilde{p}}^\varepsilon_\gamma(T -
t_k, x_k, y) d y d x_k...dx_1,
\]
where ${\widetilde{p}}^\varepsilon_\gamma(t,x,y)$ is the
fundamental solution of the parabolic equation
\[
\frac{\partial u}{\partial t} = \widetilde{L}^\varepsilon_\gamma u
:= \frac{1}{2} \Delta u +  v^\varepsilon_\gamma u
\]
and
\[
{\widetilde{Z}}^\varepsilon_{\gamma} (t,x) = \int_{ \mathbb{R}^3}
{\widetilde{p}}^\varepsilon_{\gamma}(t,x,y) d y.
\]

Since $\overline{\mathrm{P}}^x_{\gamma, T}( \omega(t_i) = 0~~{\rm
for}~{\rm some}~1 \leq i \leq k ) = 0$, in order to demonstrate
the convergence of the finite-dimensional distributions, it is
sufficient to prove that
 $ \lim_{\varepsilon \downarrow 0}
{\widetilde{Z}}^\varepsilon_{\gamma} (t,x) =
\overline{Z}_{\gamma}(t , x)$ and that for continuous functions
with compact supports $\varphi_1,...,\varphi_k$ such that ${\rm
supp} (\varphi_i) \subset \mathbb{R}^3 \setminus 0$, $1 \leq i
\leq k$, we have
\begin{equation} \label{fpu}
\lim_{\varepsilon \downarrow 0}
\int_{\mathbb{R}^3}...\int_{\mathbb{R}^3}
{\widetilde{p}}^\varepsilon_\gamma(t_1,x,x_1)...{\widetilde{p}}^\varepsilon_\gamma(t_k
- t_{k-1}, x_{k-1}, x_k)  \varphi_1(x_1)...\varphi_k(x_k)  d
x_k...dx_1 =
\end{equation}
\[
\int_{\mathbb{R}^3}...\int_{\mathbb{R}^3}
{{\overline{p}}}_\gamma(t_1,x,x_1)...{{\overline{p}}}_\gamma(t_k -
t_{k-1}, x_{k-1}, x_k)  \varphi_1(x_1)...\varphi_k(x_k)  d
x_k...dx_1.
\]
Let us start by showing that for each $k \in (0,1)$ and continuous
function $\varphi$ with compact support we have
\begin{equation} \label{lmi}
\lim_{\varepsilon \downarrow 0}||v_\varepsilon -v||_{C(F \times
[k, 1/k])} =0,
\end{equation}
where
\[
v_\varepsilon = \int_{\mathbb{R}^3}
{\widetilde{p}}^\varepsilon_\gamma(t,x,y) \varphi(y)d y,~~v =
\int_{\mathbb{R}^3} {{\overline{p}}}_\gamma(t,x,y) \varphi(y)d
y,~~F = \{x \in \mathbb{R}^3: |x| \in [k, 1/k] \}.
\]
Consider
\[
u_{\lambda,\varepsilon} = ( \mathcal{L}^\varepsilon_\gamma -
\lambda)^{-1} \varphi,~~u_\lambda = ( \mathcal{L}_\gamma -
\lambda)^{-1} \varphi,~~\lambda \notin {\rm spec}(
\mathcal{L}_\gamma).
\]
It is shown in \cite{A} that $ \mathcal{L}^\varepsilon_\gamma$
converges to $ \mathcal{L}_\gamma$ in the norm resolvent sense
when $ \lambda \notin {\rm spec}( \mathcal{L}_\gamma)$. Thus
\begin{equation} \label{str}
||u_{\lambda,\varepsilon} - u_\lambda||_{L^2( \mathbb{R}^3)}
\rightarrow 0~~{\rm as}~~\varepsilon \downarrow 0,~~\lambda \notin
{\rm spec}( \mathcal{L}_\gamma).
\end{equation}
Since there is a neighborhood of $F$ (that does not include the
origin) where
\begin{equation} \label{nmm}
(\frac{\Delta}{2} - \lambda) u_{\lambda,\varepsilon} =
(\frac{\Delta}{2} - \lambda) u_\lambda = 0,~~0 < \varepsilon \leq
\varepsilon_0,
\end{equation}
when $\varepsilon_0$ is small enough, from standard a priori
estimates for elliptic equations and (\ref{str}) it follows that
\begin{equation}
||u_{\lambda,\varepsilon} - u_\lambda||_{C(F)}  \rightarrow
0~~{\rm as}~~\varepsilon \downarrow 0,~~\lambda \notin {\rm spec}(
\mathcal{L}_\gamma).
\end{equation}

We have
\begin{equation} \label{cntr}
v_\varepsilon = \frac{1}{2 \pi i} \int_{\Gamma(a)}  (
\mathcal{L}^\varepsilon_\gamma - \lambda)^{-1} \varphi e^{\lambda
t} d \lambda,~~v = \frac{1}{2 \pi i} \int_{\Gamma(a)}  (
\mathcal{L}_\gamma - \lambda)^{-1} \varphi e^{\lambda t} d
\lambda.
\end{equation}
Since the norm of the resolvent of an operator does not exceed the
distance from the spectrum, we have
\begin{equation} \label{nmn}
||u_{\lambda,\varepsilon}||_{L^2(
\mathbb{R}^3)},~~||u_\lambda||_{L^2( \mathbb{R}^3)} \leq
{c}{|\lambda|^{-1}},~~\lambda \in \Omega_a,~~0 < \varepsilon \leq
\varepsilon_0,
\end{equation}
where $\Omega_a$ is the region in the complex $\lambda$-plane
which is to the right of the two rays starting at $\lambda = a$
and forming $\pi/4$ angles with the the ray $[a, - \infty)$. From
here it follows that the contour $\Gamma(a)$ in (\ref{cntr}) can
be replaced by $\Gamma'(a) = \partial \Omega_a$:
\begin{equation} \label{cntr2}
v_\varepsilon = \frac{1}{2 \pi i} \int_{\Gamma'(a)}
u_{\lambda,\varepsilon} e^{\lambda t} d \lambda,~~v = \frac{1}{2
\pi i} \int_{\Gamma'(a)}  u_\lambda e^{\lambda t} d \lambda.
\end{equation}
From (\ref{nmm}) and (\ref{nmn}) it follows that
\begin{equation} \label{nmn4}
||u_{\lambda,\varepsilon}||_{H^2(F)},~~||u_\lambda||_{H^2(F)} \leq
{c'},~~\lambda \in \Gamma'(a),~~0 < \varepsilon \leq
\varepsilon_0,
\end{equation}
for some constant $c'$, and the Sobolev embedding theorem implies
that
\begin{equation} \label{nmn6}
||u_{\lambda,\varepsilon}||_{C(F)},~~||u_\lambda||_{C(F)} \leq
{c''},~~\lambda \in \Gamma'(a),~~0 < \varepsilon \leq
\varepsilon_0.
\end{equation}
This allows us to split the contour of integration in the
integrals in (\ref{cntr2}) into two parts: a part where $|\lambda|
>> 1$ and a bounded contour. The contribution to $v_\varepsilon$
and $v$ from the first part can be made arbitrarily small. The
Lebesgue dominated convergence theorem implies that the difference
of the integrals over the bounded part of the contour tends to
zero. Thus we have (\ref{lmi}).

It is easy to see that (\ref{lmi}) implies (\ref{fpu}). Let us now
prove that $ \lim_{\varepsilon \downarrow 0}
{\widetilde{Z}}^\varepsilon_{\gamma} (t,x) =
\overline{Z}_{\gamma}(t , x)$. Let $\xi$ be an infinitely smooth
function that is equal to zero for $|x| \leq 1$ and equal to one
for $|x| \geq 2$. Let
\[
{\widetilde{Z}}^\varepsilon_{\gamma} =
{\widetilde{Z}}^{\varepsilon,1}_{\gamma} +
{\widetilde{Z}}^{\varepsilon,2}_{\gamma},~~\overline{Z}_{\gamma} =
\overline{Z}^1_{\gamma}+\overline{Z}^2_{\gamma},
\]
where
\[
{\widetilde{Z}}^{\varepsilon,1}_{\gamma} (t,x) = \int_{
\mathbb{R}^3}  {\widetilde{p}}^\varepsilon_{\gamma}(t,x,y) \xi(y)
d y,~~{\widetilde{Z}}^{\varepsilon,2}_{\gamma} (t,x) = \int_{
\mathbb{R}^3} {\widetilde{p}}^\varepsilon_{\gamma}(t,x,y) (1-
\xi(y)) d y,
\]
and $ \overline{Z}^1_{\gamma}$, $\overline{Z}^2_{\gamma}$ are
defined similarly. Then $ \lim_{\varepsilon \downarrow 0}
{\widetilde{Z}}^{\varepsilon,2}_{\gamma} (t,x) =
\overline{Z}_{\gamma}^2(t , x)$ by (\ref{lmi}) with $\varphi  = 1
- \xi$.

In order to estimate $|\widetilde{Z}^{\varepsilon,1}_{\gamma}
(t,x) - \overline{Z}_{\gamma}^1(t , x)|$, we consider the
following problem:
\begin{equation} \label{g5}
\frac{\partial u}{\partial t} - \mathcal{L}^\varepsilon_\gamma u =
h,~~u(0,x) = \xi(x),
\end{equation}
where $h \in C_0^{\infty}( \mathbb{R}^3)$ depends on $x$ only. The
solution of this equation has the form
\[
u = {\widetilde{Z}}^{\varepsilon,1}_{\gamma} - \frac{1}{2 \pi i}
\int_{\Gamma'(a)} \lambda^{-1} (\mathcal{L}^\varepsilon_\gamma -
\lambda)^{-1} h e^{\lambda t} d \lambda.
\]
Note that $u(t,x) = \xi(x)$ satisfies (\ref{g5}) with $h = -\Delta
\xi /2 $. Thus
\[
\xi = {\widetilde{Z}}^{\varepsilon,1}_{\gamma} + \frac{1}{4 \pi i}
\int_{\Gamma'(a)} \lambda^{-1} (\mathcal{L}^\varepsilon_\gamma -
\lambda)^{-1} (\Delta \xi )  e^{\lambda t} d \lambda.
\]
The same formula is valid with $ \overline{Z}_{\gamma}^1$ instead
of $ \widetilde{Z}^{\varepsilon,1}_{\gamma}$ and
$\mathcal{L}_\gamma$ instead of $\mathcal{L}^\varepsilon_\gamma$.
Therefore,
\[
|\widetilde{Z}^{\varepsilon,1}_{\gamma} (t,x) -
\overline{Z}_{\gamma}^1(t , x)| =\frac{1}{4 \pi}|
\int_{\Gamma'(a)} \lambda^{-1} (\mathcal{L}^\varepsilon_\gamma -
\lambda)^{-1} (\Delta \xi )  e^{\lambda t} d \lambda  -
\int_{\Gamma'(a)} \lambda^{-1} (\mathcal{L}_\gamma - \lambda)^{-1}
(\Delta \xi )  e^{\lambda t} d \lambda|.
\]
Note that $\Delta \xi$ has compact support, and therefore we can
treat the right hand side in this formula in the same way as we
treated the difference between $v_\varepsilon$ and $v$ thus
obtaining  that $ \lim_{\varepsilon \downarrow 0}
{\widetilde{Z}}^{\varepsilon,1}_{\gamma} (t,x) =
\overline{Z}^1_{\gamma}(t , x)$.

It remains to note that the family of measures
${{\mathrm{P}}}_{\gamma, T}^{x,\varepsilon}$, $0< \varepsilon \leq
1$, is tight, as follows from Lemma~\ref{lemtight}. \qed
\\
\\
{\bf Remark.} It follows from the arguments above that for any
$c,t
> 0$,
 there is a constant $C$ such that ${\widetilde{Z}}^\varepsilon_{\gamma}
 (t,x) \leq C$ provided that $\varepsilon \leq c$,  $1/c \leq |x| \leq c$ and $\gamma \leq
 c$.
 \\

Next we prove two statements that are needed for the proofs of
Theorems~\ref{convlem} and \ref{pro1}. Note that we already
demonstrated the convergence of partition functions and
finite-dimensional distributions in Theorem~\ref{convlem}, and
using these facts in the proofs of Lemmas~\ref{lemtight} and
\ref{lemtight2} does not involve a circular argument.
\begin{lemma} \label{lemtight}
Let $v_\gamma^\varepsilon$ be given by (\ref{poten}), the
corresponding Hamiltonian $H^{\varepsilon}_{\gamma,T}$ given by
(\ref{hamt}) with $v_{\gamma}^\varepsilon$ instead of $v$, and the
Gibbs measure ${{\mathrm{P}}}_{\gamma, T}^{x,\varepsilon}$ given
by~(\ref{gibm}) with $\beta = 1$ and $H^{\varepsilon}_{\gamma,T}$
instead of $H_T$. For each $T >0$ and $c \in \mathbb{R}$, the
family of measures $\{ {{\mathrm{P}}}_{\gamma, T}^{x,\varepsilon},
\varepsilon \leq c, \gamma \leq c, |x| \leq c \}$ is tight.
\end{lemma}
\proof
 To prove tightness it
is enough to demonstrate that for each $\eta, \alpha > 0$ there is
$0 < \delta < 1$ such that for all $u \in [0,T]$ we have
\begin{equation} \label{crite}
\mathrm{P}^{x,\varepsilon}_{\gamma, T} (\sup_{u \leq s \leq
\min(u+\delta, T)} |\omega(s) - \omega(u)| > \alpha) \leq \delta
\eta.
\end{equation}
Let $\eta, \alpha >0$ be fixed. If $\varepsilon_0 > 0$, the
density of the measure $\mathrm{P}^{x,\varepsilon}_{\gamma, T}$
with respect to the measure induced by the Brownian motion
starting at $x$ is bounded from above uniformly in $\varepsilon_0
\leq \varepsilon \leq c$, $\gamma \leq c$, $|x| \leq c$. Thus we
can find $0 < \delta < 1$ such that (\ref{crite}) holds for
$\varepsilon_0 \leq \varepsilon \leq c$. Therefore, it is
sufficient to prove (\ref{crite}) under the assumption that
$\varepsilon$ is sufficiently small so that $
v_\gamma^\varepsilon(x) = 0$ for $|x| \geq \alpha/8$.

Consider the following events in $C([0,T], \mathbb{R}^3)$,
\[
\widetilde{\mathcal{E}}_\delta^T = \{ \omega: \sup_{u \leq s \leq
\min(u+\delta, T)} |\omega(s) - \omega(u)| \geq
\alpha\},~~~\mathcal{E}_\delta^T = \{\omega:  \sup_{s \leq
\min(\delta,T)} |\omega(s) - \omega(0)| \geq \alpha/8 \}.
\]
For a continuous function $\omega: [0,T] \rightarrow
\mathbb{R}^3$, let
\[
\tau = \min(T, \inf\{t \geq 0: |\omega(t)| \geq \alpha/4 \}).
\]
Then
\[
\mathrm{P}^{x,\varepsilon}_{\gamma, T}(
\widetilde{\mathcal{E}}_\delta^T) =
(\widetilde{Z}^\varepsilon_{\gamma}(T,x))^{-1}
\mathrm{E}^x_{T}(\exp(\int_0^{T} v_\gamma^\varepsilon (\omega(t))
d t) \chi_{ \widetilde{\mathcal{E}}_\delta^T}) \leq
\]
\[
(\widetilde{Z}^\varepsilon_{\gamma}(T,x))^{-1}
\mathrm{E}^x_{T}\left(\exp(\int_0^{\tau} v_\gamma^\varepsilon
(\omega(t)) d t) \mathrm{E}^{\omega(\tau)}_{T-\tau} ( \chi_{
\mathcal{E}^\tau_\delta} \exp( \int_0^{T -\tau}
v_\gamma^\varepsilon (\omega(t)) d t ) ) \right),
\]
where $ \mathrm{E}^{x}_{T}$ denotes the expectation with respect
to the measure induced by the Brownian motion starting at the
point $x$.
 Since
\[
\mathrm{E}^x_{T} \exp(\int_0^{\tau} v_\gamma^\varepsilon
(\omega(t)) d t)  \leq \widetilde{Z}^\varepsilon_{\gamma}(T,x)
\]
and
\[
 \mathrm{E}^{\omega(\tau)}_{T-\tau} ( \chi_{
\mathcal{E}^\tau_\delta} \exp( \int_0^{T -\tau}
v_\gamma^\varepsilon (\omega(t)) d t ) ) \leq \sup_{ \alpha/4 \leq
|x| \leq c} \mathrm{E}^{x}_{T} ( \chi_{ \mathcal{E}^T_\delta}
\exp( \int_0^{T} v_\gamma^\varepsilon (\omega(t)) d t ) ),
\]
we have
\[
\mathrm{P}^{x,\varepsilon}_{\gamma, T}(
\widetilde{\mathcal{E}}_\delta^T)  \leq  \sup_{ \alpha/4 \leq |x|
\leq c} \mathrm{E}^{x}_{T} ( \chi_{ \mathcal{E}^T_\delta} \exp(
\int_0^{T} v_\gamma^\varepsilon (\omega(t)) d t ) ).
\]
 Let
\[ \sigma = \min(\delta, \inf\{t \geq 0: |\omega(t) - \omega(0)| = \alpha/8\}
).
\]
Then, since $ v_\gamma^\varepsilon(x) = 0$ for $|x| \geq
\alpha/8$,
\[
 \sup_{\alpha/4 \leq |x| \leq c} \mathrm{E}^{x}_{T} ( \chi_{
\mathcal{E}^T_\delta} \exp( \int_0^{T} v_\gamma^\varepsilon
(\omega(t)) d t ) ) \leq
 \sup_{\alpha/4 \leq |x| \leq c} \mathrm{E}^{x}_{T} ( \chi_{
\mathcal{E}^T_\delta} \mathrm{E}^{\omega(\sigma)}_{T - \sigma}
\exp( \int_0^{T-\sigma} v_\gamma^\varepsilon (\omega(t)) d t ) ).
\]
Note  that
\[
\mathrm{E}^{\omega(\sigma)}_{T - \sigma} \exp( \int_0^{T-\sigma}
v_\gamma^\varepsilon (\omega(t)) d t ) \leq \sup_{ \alpha/8 \leq
|x| \leq c+ \alpha/8 } \mathrm{E}^{x}_{T} \exp( \int_0^{T}
v_\gamma^\varepsilon (\omega(t)) d t ) =
\]
\[
\sup_{ \alpha/8 \leq |x| \leq c+ \alpha/8 }
\widetilde{Z}^\varepsilon_{\gamma}(T,x) \leq
 c(\alpha)
\]
for some constant $c(\alpha)$, where the last inequality is due to
the Remark following the proof of Theorem~\ref{convlem}.
 It then remains to choose $\delta$
such that $\mathrm{E}^{x}_{T} ( \chi_{ \mathcal{E}^T_\delta} )
\leq \delta \eta /c(\alpha)$, and estimate  (\ref{crite}) follows.
\qed

\begin{lemma} \label{lemtight2}
Suppose that $\gamma(T)$ is bounded and such that $\gamma(T)
\sqrt{T} \rightarrow +\infty$ as $T \rightarrow +\infty$. Suppose
that $s > 0$, $S(T)$ satisfies (\ref{sttt}) and $T_0$ is such that
$S(T) + s/\gamma^2(T) \leq T$ for $T \geq T_0$. Let the process
$y^T(t)$, $0 \leq t \leq s$, be defined on the probability space
$C([0,T], \mathbb{R}^3)$ with the measure
$\overline{\mathrm{P}}^x_{\gamma(T), T}$ by $y^T(t) = \gamma(T)
\omega(S(T)+t/\gamma^2(T))$, $0 \leq t \leq s$. Then, for each $K
> 0$,  the family of measures on $C([0,s], \mathbb{R}^3)$ induced by
the processes $y^T$ with $T \geq T_0$ and $0 < |x| \leq K$ is
tight.
\end{lemma}
\proof Due to self-similarity, it is sufficient to consider the
case when $\gamma(T) = 1$.
%
First observe that for any $\eta > 0$ there is $a > 0$ such that
\[
\overline{\mathrm{P}}^x_{1, T}(|y^T(0)| > a)  \leq \eta
\]
for $T \geq T_0$ and $0 < |x| \leq K$. This easily follows from
Lemma~\ref{bnbn1} in the same way that the convergence of
finite-dimensional distributions was demonstrated in
Theorem~\ref{pro1}.

For a continuous function $\omega: [0,T] \rightarrow
\mathbb{R}^3$, let
\[
m^T(\omega, \delta) = \sup_{|t_1 - t_2| \leq \delta,~S(T) \leq
t_1, t_2 \leq S(T) + s} |\omega(t_1) - \omega(t_2)|.
\]
To prove tightness, it is sufficient to show that for each
$\alpha, \eta >0$, there is $\delta > 0$ such that
\begin{equation} \label{tpr}
\overline{\mathrm{P}}^x_{1, T}(m^T(\omega, \delta) > \alpha) \leq
\eta
\end{equation}
 for $T \geq T_0$ and $0 < |x| \leq K$. By
Theorem~\ref{convlem},
\begin{equation} \label{hh1}
\overline{\mathrm{P}}^x_{1, T}(m^T(\omega, \delta) > \alpha) \leq
\liminf_{\varepsilon \downarrow 0} \mathrm{P}^{x,\varepsilon}_{1,
T}(m^T(\omega, \delta) > \alpha).
\end{equation}
For a continuous function $\omega: [0,T] \rightarrow
\mathbb{R}^3$, let
\[
\tau = \min(S(T) + s, \inf\{t \geq S(T): |\omega(t)| \geq \alpha/2
\}).
\]
For $0 \leq a \leq b \leq r$, let
\[
\mathcal{E}^r_{a,b} = \{ \omega \in C([0,r], \mathbb{R}^3) :
\sup_{|t_1 - t_2| \leq \delta,~a \leq t_1, t_2 \leq b}
|\omega(t_1) - \omega(t_2)| > \alpha \}.
\]
 Then
\[
\liminf_{\varepsilon \downarrow 0} \mathrm{P}^{x,\varepsilon}_{1,
T}(m^T(\omega, \delta) > \alpha) =
(\widetilde{Z}^\varepsilon_{\gamma}(T,x))^{-1}
\mathrm{E}^x_{T}(\exp(\int_0^{T} v_\gamma^\varepsilon (\omega(t))
d t) \chi_{\mathcal{E}^T_{S(T),S(T)+s}}) \leq
\]
\begin{equation} \label{hh2}
(\widetilde{Z}^\varepsilon_{\gamma}(T,x))^{-1}
\mathrm{E}^x_{T}\left(\exp(\int_0^{\tau} v_\gamma^\varepsilon
(\omega(t)) d t) \mathrm{E}^{\omega(\tau)}_{T - \tau}
(\exp(\int_0^{T-\tau} v_\gamma^\varepsilon (\omega(t)) d t)
\chi_{\mathcal{E}^{T-\tau}_{0,S(T)+s -\tau}})\right) \leq
\end{equation}
\[
(\widetilde{Z}^\varepsilon_{\gamma}(T,x))^{-1}
(\widetilde{Z}^\varepsilon_{\gamma}(S(T) + s,x)) \sup_{|x| \geq
\alpha/2} \mathrm{E}^{x}_{T - S(T)} (\exp(\int_0^{T-S(T) }
v_\gamma^\varepsilon (\omega(t)) d t)
\chi_{\mathcal{E}^{T-S(T)}_{0,s}}).
\]
Since $v_\gamma^\varepsilon (x) = 0$ when $x \geq \alpha/4$,
provided that $\varepsilon$ is small enough, we have by the Markov
property
\begin{equation} \label{hh3}
\sup_{|x| \geq \alpha/2} \mathrm{E}^{x}_{T - S(T)}
(\exp(\int_0^{T-S(T) } v_\gamma^\varepsilon (\omega(t)) d t)
\chi_{\mathcal{E}^{T-S(T)}_{0,s}}) \leq
\end{equation}
\[
\sup_{|x| \geq \alpha/4}
\widetilde{Z}^\varepsilon_{\gamma}(T-S(T),x)
\mathrm{P}^0_s(\{\omega: \sup_{|t_1 - t_2| \leq \delta,~0 \leq
t_1, t_2 \leq s} |\omega(t_1) - \omega(t_2)| > \alpha  \}).
\]
Recall from the proof of Theorem~\ref{convlem} that
 $ \lim_{\varepsilon \downarrow 0}
{\widetilde{Z}}^\varepsilon_{\gamma} (t,x) =
\overline{Z}_{\gamma}(t , x)$. Therefore, combining (\ref{hh1}),
(\ref{hh2}) and (\ref{hh3}) we see that
\[
\overline{\mathrm{P}}^x_{1, T}(m^T(\omega, \delta) > \alpha) \leq
(\overline{Z}_{\gamma}(T,x))^{-1} (\overline{Z}_{\gamma}(S(T) +
s,x)) \times
\]
\[
 \sup_{|x| \geq
\alpha/4} \overline{Z}_{\gamma}(T-S(T),x) \mathrm{P}^0_s(\{\omega:
\sup_{|t_1 - t_2| \leq \delta,~0 \leq t_1, t_2 \leq s}
|\omega(t_1) - \omega(t_2)| > \alpha  \}).
\]
Formula (\ref{tpr}) now follows from Lemma~\ref{bnbn1}. \qed
\\
\\

\noindent {\bf \large Acknowledgements}: While working on this
article, M. Cranston was supported by NSF grant  DMS-0706198, L.
Koralov was supported by NSF grant DMS-0706974, S. Molchanov and
B. Vainberg were supported by NSF grant DMS-0706928.


\begin{thebibliography}{9}


\bibitem{A} S. Albeverio, F. Gesztesy, R. H{\o}egh-Krohn, H.
              Holden.
\newblock{\em Solvable models in quantum mechanics},
 AMS Chelsea Publishing, Providence, RI, 2005.

\bibitem{AS}
K. Alexander, V. Sidoravicius.
 \newblock{\em Pinning of polymers and interfaces by random potentials},
Ann. Appl. Probab., 16, 2006, 2, 636--669.

\bibitem{BS} E. Bolthausen, U. Schmock. {\it On self-attracting
$d$-dimensional random walks}, Ann. Probab. Volume 25, Number 2
(1997), 531-572.

\bibitem{CH}
P. Carmona,  Y. Hu.
\newblock{\em On the partition function of a directed polymer in a
              {G}aussian random environment},
 Probab. Theory Related Fields, 124, 2002, 3, 431--457.

\bibitem{CGH}
P. Carmona, F. Guerra,  Y. Hu, O. Menjane.
\newblock{\em Strong disorder for a certain class of directed polymers in a
random environment}, J. Theoret. Probab., 19, 2006, 1, 134--151.

\bibitem{CSY1}
F. Comets, T. Shiga, N. Yoshida. \newblock{\em Directed polymers
in a random environment: path localization and strong disorder},
Bernoulli. Official Journal of the Bernoulli Society for
Mathematical Statistics and Probability, 9, 2003, 4, 705--723.

\bibitem{CSY2}
F. Comets, T. Shiga, N. Yoshida.\newblock{\em Probabilistic
analysis of directed polymers in a random environment: a review},
Stochastic analysis on large scale interacting systems, Adv. Stud.
Pure Math., 39, 115--142, Math. Soc. Japan, Tokyo, 2004.

\bibitem{CY} F. Comets, N. Yoshida.
\newblock{\em Brownian directed polymers in random environment},
Comm. Math. Phys., 254, 2005, 2, 257--287.

\bibitem{CKMV}
M. Cranston, L. Koralov, S. Molchanov, B. Vainberg.
\newblock{\em A continuous model for homopolymers},
to appear in Journal of Functional Analysis.


\bibitem{CM}
M. Cranston, S. Molchanov.
\newblock{\em Analysis of  homopolymer,} submitted.

\bibitem{D}
R. L. Dobrushin, S. Kusuoka.
\newblock{\em Statistical mechanics and fractals},
Lecture Notes in Mathematics, 1567, Springer-Verlag, Berlin, 1993.

\bibitem{G}
G. Giacomin. \newblock{\em Random polymer models}, Imperial
College Press, London, 2007.

\bibitem{K}
H. Kesten.
\newblock{\em Scaling relations for {$2$}{D}-percolation},
Comm. Math. Phys., 109, 1987, 1, 109--156.

\bibitem{KW} W. Konig. {\it Self-repellent and self-attractive
path measures in statistical mechanics}, Habilitationsschrift,
Techniche Universitat Berlin, 2000.

\bibitem{Ko} L. Koralov. {\it On rotation-invariant diffusion
processes}, in preparation.

\bibitem{LL} L. D. Landau, E.M. Lifshitz.
\newblock{\em Statistical Physics (Course of Theoretical Physics, Volume 5},
Oxford, Boston, 2000.

\bibitem{LGK} I. M. Lifshitz,   A. Yu. Grosberg,   A. R. Khokhlov.
\newblock{\em Some problems of the statistical physics of polymer chains
with volume interaction}, Rev. Modern Phys.,50,1978, 3,683-713.











\end{thebibliography}
\end{document}